
\input amstex
\documentstyle{amsppt}
\magnification=\magstep1
\NoBlackBoxes


\topmatter

\title Inverse formula for the Blaschke-Levy 
representation with applications to zonoids 
and sections of star bodies.
\endtitle

\author Alexander Koldobsky \endauthor 
\address Division of Mathematics and Statistics, 
University
of Texas at San Antonio, San Antonio, TX 78249, U.S.A. \endaddress
\email koldobsk\@ringer.cs.utsa.edu \endemail

\abstract  We say that an even continuous function $H$ 
on the unit sphere $\Omega$ 
in $R^n$ admits the Blaschke-Levy representation with $q>0$ 
if there exists an even function $b\in L_1(\Omega)$ so that
$H^q(x)=\int_\Omega  |(x,\xi)|^q b(\xi)\ d\xi$ 
for every $x\in \Omega.$ This representation has numerous 
applications in convex geometry, probability and Banach space 
theory. In this paper, we present a simple formula 
(in terms of the derivatives of $H$)
for calculating $b$ out of $H.$ We use this formula 
to give a sufficient condition for isometric embedding of a space 
into $L_p$ which 
contributes to the 1937 P.Levy's problem and to the study of zonoids. 
Another application gives a Fourier transform formula for the volume of 
$(n-1)$-dimensional central sections of star bodies in $R^n.$ We apply 
this formula to find the minimal and maximal volume of central 
sections of the unit balls of the spaces $\ell_p^n$ with $0<p<2.$ 
\endabstract


\rightheadtext{Inverse formula for the Blaschke-Levy representation}

\endtopmatter \document \baselineskip=14pt

\head 1. Introduction \endhead

For $q>0,$ we say that an even continuous function $H$
on $\Bbb R^n$  admits the Blaschke-Levy representation 
with the exponent $q$ if there exists an even function 
$b$ on the unit sphere $\Omega$ in $\Bbb R^n$ so that
$b\in L_1(\Omega)$ and, for every $x\in \Bbb R^n,$
$$H(x) = \int_\Omega |(x,\xi)|^q b(\xi)\ d\xi,\tag{1}$$
where $(x,\xi)$ stands for the scalar product.

It was known to Blaschke \cite{3} that every infinitely
differentiable function on the sphere admits the representation
(1) with $q=1.$ On the other hand, the representation 
(1) is known in the probability theory under the name of P.Levy,
and it was an important part of P.Levy's theory of stable 
processes \cite{19} that the function $\|x\|^q$ admits the representation 
(1) with a measure in place of the function $b,$ 
where $(\Bbb R^n,\|\cdot\|)$ is any $n$-dimensional subspace of $L_q.$ 
In mathematical physics the representation (1) is called
the plain-wave expansion.
\bigbreak
The Blaschke-Levy representation has had numerous applications
to convex geometry, probability and Banach space theory.
One of the most popular ways to apply the Blaschke-Levy
representation is based on the fact that the representation
is unique for every $q>0$ which is not an even integer
(the uniqueness fails if $q$ is an even integer, because 
only a finite number of moments of the functions must be equal).
The uniqueness was first shown by Blaschke \cite{3} 
in the case where $q=1$ and $n=3.$ Aleksandrov
\cite{1} proved the uniqueness for $q=1$ and arbitrary dimension, 
and P. Levy \cite{19} did it for $0<q<2.$
The last two results are valid for signed measures in place of $b.$
The uniqueness for every $q$ which is not an even integer was
established by Kanter \cite{13}. For different proves and 
applications of the uniqueness theorem see \cite{11, 20, 23, 24, 14}.  
In Section 2 we present a Fourier transform proof which 
is close to that from \cite{14}.

\bigbreak

The existence of the representation (1) with $q=1$ for infinitely 
differentiable functions was known to Blaschke \cite{3}. A precise 
proof under a weaker assumption that $H\in C^{n+2}(\Omega)$
was given by Schneider \cite{26} who found a spherical harmonics
expansion for the function $b$ (which turned out to be a continuous 
function on $\Omega.)$ Later Goodey and Weil \cite{10}
proved the existence of the representation (1) (also with $q=1)$ 
for the functions $H$ of the class 
$C^{(n+5)/2}$ where the function $b$ appears to belong to the space 
$L_2(\Omega).$ Weil \cite{28} found a generating distribution  
for the support function of any centered convex body.
Richards \cite{25} showed that the representation (1) exists for any
$q\in (0,2)$ and any $H\in C^{n+q+1}(\Omega).$ A generalization 
of this result to the case of arbitrary $q>0$ which is not an 
even integer was given in \cite{16}. All the results mentioned above
were based on the use of spherical harmonics.
A connection between the Blaschke-Levy representation and the
Fourier transform was found in \cite{14} where it was shown that
the function $b$ is the restriction to the sphere $\Omega$ of the 
Fourier transform of $H$ (we present
a short version of that proof in Theorem 1 below; in fact,
in \cite{14} the  Fourier transform of $H$ was restricted 
to a hyperplane). This fact was used to show that every
norm in $\Bbb R^n$ admits the Blaschke-Levy representation with
every $q>0$ which is not an even integer, but we must allow $b$
to be a distribution and 
the representation (1) is considered in a generalized form. 
Note that the Fourier transform 
connection was used in \cite{14, 15} to obtain exact representations
for certain norms, which, in particular, led to
applications to positive definite functions and embedding of Banach
spaces.

\bigbreak

A remarkable feature of Schneider's spherical harmonics
construction is that it allows to gain control over the
function $b$ by estimating the $L_\infty(\Omega)$-norm
of $b$ in terms of $H.$ Namely,
Schneider \cite{26} showed that, for any $H\in C^{n+2}(\Omega),$
the function $b$ appearing in the Blaschke-Levy representation 
with $q=1$ satisfies the inequality
$$\|b\|_\infty \leq K\|H\|_{L_{2}(\Omega)} + 
L\|\Delta_\Omega^{2r}H\|_{L_{2}(\Omega)},\tag{2}$$
where $\Delta_\Omega$ is the Laplace-Beltrami operator,
$r>(n+2)/2$ and $K$ and $L$ are constants which are
given as the sum of certain series'. Schneider \cite{27}
used this inequality to construct non-trivial zonoids 
whose polars are zonoids. 
In order to do that, he considered
a perturbation of the Euclidean norm by means of an infinitely
differentiable function $f$ on the sphere $\Omega:$ 
put $H(x)=\|x\|_2 (1+\lambda f(x/\|x\|_2)),\ x\in \Bbb R^n,$ where
$\|\cdot\|_2$ is the Euclidean norm and $\lambda$ is a (small)
real number. If the function $b$ corresponding to $H$
in the representation (1) is non-negative,
then $H$
is the norm of a subspace of $L_1,$ and, therefore,
it is the support function of a body whose polar is a zonoid.
Since the function $b$ corresponding to the Euclidean norm
$\|x\|_2$ in the Blaschke-Levy representation is a constant,
and the $\ell_\infty$-norm of the perturbing function $b$ is 
controlled by $\lambda$ because of (2), 
one can choose $\lambda$ small enough so that the function
$b$ corresponding to $H$ is non-negative. It is easy to see that
making $\lambda$ even smaller (if necessary) one can make the body
$\{x: H(x)\le 1\}$ to be a zonoid too.

The inequality (2) was generalized in \cite{16} to the case
of the Blaschke-Levy representation with any $q>0$ which
is not an even integer. This led to a construction of
common subspaces of $L_q$-spaces: for any $n\in N$
and any compact subset $Q$ of $(0,\infty)\setminus
\{even\ integers\},$ 
there exists an $n$-dimensional non-Hilbertian Banach space 
which is isometric to a subspace of $L_q$ for every $q\in Q.$

\bigbreak

This paper is an attempt to gain more control over the function $b$
by presenting an inverse formula for the representation (1) which
does not involve spherical harmonics or the Fourier transform, and by
giving a simpler version of the inequality (2) with computable 
constants. We start with the Fourier transform inverse formula
showing that $b$ is the restriction to the sphere of the Fourier
transform of the function $H$ (which is homogeneous 
of degree $q$ because of (1)). However, to avoid the calculation of the
Fourier transform, we first apply the Laplace operator to the
function $H$ as many times as it is necessary to make the result
homogeneous of degree less or equal than $-n+1.$ Note that action of 
the Laplace operator does not change the restriction of the Fourier
transform to the sphere (up to a sign). The crucial point is that,
by Lemmas 3 and 4, the Fourier transform of a homogeneous function
of degree less or equal than $-n+1$ can easily be expressed in
terms of the function itself. 

In this way we show that, for 
every $q>0$ which is not an integer and every even homogeneous
function $H$ of degree $q$ on $\Bbb R^n$ such that the restriction
to sphere $H\vert_\Omega$ belongs to the space $C^{n+[q]}(\Omega),$
there exists the Blaschke-Levy representation with the exponent $q,$
where the corresponding function $b$ is given by
$$b(\xi) = (-1)^k {\pi\over{2(2\pi)^{n-1}C_{-n-q+2k}C_q}} 
\int_\Omega |(\theta,\xi)|^{-n-q+2k} 
(\Delta^k H^q)(\theta)\ d\theta,$$
for every $\xi\in \Omega,$   
where $k=(n+[q])/2$ if $n+[q]$ is an even integer, and  
$k=(n+[q]+1)/2$ if $n+[q]$ is an odd integer. 

If $q$ is an odd integer and the dimension $n$ is an even integer
the expression for $b$ is as follows:
$$b(\xi) = {{(-1)^{(n+q-1)/2}\pi}\over{(2\pi)^{n-1}C_q }}
\int_{\Omega\cap \{(\theta,\xi)=0\}} 
\Delta^{(n+q-1)/2} H^q(\theta)\ d\theta.$$

If both $q$ and $n$ are odd integers the technique of this paper
does not work for the reason that, in this case, the Laplace 
transform of $H$ may contain a part supported at zero.
As it was mentioned above, if $q$ is an even integer the 
uniqueness fails.

\bigbreak

In Section 4 we apply the inverse formulae to get
a new criterion for the existence of
an isometric embedding of a given space into $L_q.$ 
Finding such criteria is a matter of the 1937 P. Levy's
problem (see \cite{19}).
We calculate
the functions $b$ for certain perturbations of the
Euclidean norm, and show the way  to get exact constants 
$\lambda$ in Schneider's construction.

\medbreak

In Section 5 we use our results to get a Fourier transform
formula for the volume of central $(n-1)$-dimensional sections
of centrally symmetric star bodies in $\Bbb R^n.$ If $K$ is such a body
then, for every $\xi\in \Omega,$
$${Vol}_{n-1}(K\cap \xi^{\bot}) = {1\over{\pi (n-1)}}
(\|x\|^{-n+1})^{\wedge}(\xi)$$
where $\xi^{\bot}= \{x\in \Bbb R^n : (x,\xi)=0\},$
and $\|x\|=\min\{a\in \Bbb R; ax\in K\}.$
Finally, we use this formula to show that
the minimal volume of central sections of the unit ball 
of the space $\ell_p^n,\ p\in (0,2)$ occurs if the section
is perpendicular to the vector $\xi=(1,1..,1).$  
This result proves a conjecture of Meyer and Pajor \cite{21}.

\bigbreak

\head 2. Connection between the Blaschke-Levy representation 
and the Fourier transform. \endhead 

The main tool of this paper is the Fourier transform 
of distributions. As usual, we denote by $\Cal S(\Bbb R^n)$ 
the space of rapidly decreasing infinitely differentiable 
functions (test functions) in $\Bbb R^n,$ and 
$\Cal S^{'}(\Bbb R^n)$ is
the space of distributions over $\Cal S(\Bbb R^n).$  
The Fourier transform of a distribution $f\in \Cal S^{'}(\Bbb R^n)$ 
is defined by $(\hat{f},\hat{\phi}) = (2\pi)^n (f,\phi)$ for every
test function $\phi.$  A distribution is called even homogeneous
of degree $q\in \Bbb R$  if 
$\big(f(x), \phi(x/\alpha)\big) = |\alpha|^{n+q} (f,\phi)$
for every test function $\phi$ and every 
$\alpha\in \Bbb R,\ \alpha\neq 0.$ The Fourier transform of 
an even homogeneous distribution of degree $q$ is an even 
homogeneous distribution of degree $-n-q.$

If $q>-1$ and $q$ is not an even integer, then 
the Fourier transform of the function
$h(z)=|z|^q,\ z\in \Bbb R$ is equal to 
$(|z|^q)^{\wedge}(t) = C_q |t|^{-1-q}$ 
(see \cite{8, p. 173}), where
$$C_q = {{2^{q+1}\sqrt{\pi}\ \Gamma((q+1)/2)}\over{\Gamma(-q/2)}}.$$

\medbreak

Throughout the paper, we use the following fact which is
a simple consequence of the connection
between the Fourier transform and the Radon transform.

\proclaim{Lemma 1} 
Let $q>-1,$ $q$ is not an even integer.
Then for every even test function $\phi$ with $0\notin supp(\phi)$
and every fixed vector $\xi\in\Bbb R^n$, $\xi\not = 0$, we have
$$
\int_{\Bbb R^n} |(x,\xi)|^q\ \hat{\phi}(x)\ dx=
(2\pi)^{n-1}C_q \int_{\Bbb R} |t|^{-1-q} \phi(t\xi)\ dt.\tag{3}$$
\endproclaim

\demo{Proof} 
By the well-known connection between the Fourier transform
and the Radon transform (see \cite{12}),
the function $t\to (2\pi)^n \phi(-t\xi)$ 
is the Fourier transform of the function 
$z\to \int_{(x, \xi)=z} \hat\phi (x)\, dx.$ 
(Recall that $(\hat{\phi})^{\wedge}(x) = (2\pi)^n \phi(-x).)$
Using this fact and the Fubini theorem,
for every test function $\phi$ with $0\notin supp(\phi),$
we get
$$
\int_{\Bbb R^n} |(x,\xi)|^q\ \hat\phi(x)\, dx=
\int_{\Bbb R} |z|^q
\Big(\int_{(x,\xi)=z} \hat \phi (x)\, dx\Big)\, dz=
\Big( |z|^q, \int_{(x, \xi) =z} \hat \phi (x)\, dx\Big)=$$
$${1\over{2\pi}} \Big( C_q|t|^{-1-q}, (2\pi)^n \phi(-t\xi)\Big)=
(2\pi)^{n-1}C_q \int_{\Bbb R} |t|^{-1-q} \phi(t\xi)\ dt.\qed$$ 
\enddemo

\bigbreak

\subheading{Remark 1} If $q>-1$ and $\mu$ is a Borel  
signed measure with bounded variation on $\Omega,$ then the integral 
$G(x)=\int_\Omega |(x, \xi)|^q\ d\mu(\xi)$
converges for almost all $x\in \Omega$ with respect to the
uniform measure on $\Omega.$ This follows from the fact that,
for $q>-1$ and any $\xi \in \Omega,$ 
$$W_q = \int_\Omega |(x,\theta)|^q\ dx =
2\Gamma((q+1)/2)\pi^{(n-1)/2}/\Gamma((n+q)/2) < \infty,$$
and, therefore, the restriction $G\vert_\Omega$ of the function
$G$ to $\Omega$ satisfies
$$\|G\vert_\Omega\|_1 \le \int_\Omega \ d|\mu|(\xi)
\int_\Omega |(x,\theta)|^q\ dx = W_q |\mu|(\Omega).$$
If $\mu$ has the density $f\in L_\infty(\Omega)$ then
$\|G\vert_\Omega\|_\infty \le W_q \|f\|_\infty.$
We denote by $\|\cdot\|_1$ and  $\|\cdot\|_\infty$  the
norms of the spaces $L_1(\Omega)$ and $L_\infty(\Omega),$ 
respectively.

\bigbreak

Let us calculate the Fourier transform of the function $G$
from Remark 1.

\proclaim{Lemma 2} Let $q>-1,$ $q$ is not an even integer,
and let $\mu$ be a Borel symmetric signed measure with bounded variation 
on $\Omega.$ Then the Fourier transform $\hat{G}$ of the function
$G(x)= \int_\Omega |(x,\xi)|^q \ d\mu(\xi)$ has the property that
for every even test function $\phi$ with $0\notin supp(\phi),$ 
$$(\hat{G},\phi) = (2\pi)^{n-1} C_q
\int_\Omega d\mu(\xi) \int_{\Bbb R}|t|^{-1-q} \phi(t\xi)\ dt. \tag{4}$$
\endproclaim

\demo{Proof} By Remark 1, $G$ is an even homogeneous function of 
degree $q$ whose restriction to the sphere belongs to the space 
$L_1(\Omega).$ 
For every even test function $\phi$ with $0\notin supp(\phi),$
using Lemma 1  and the Fubini theorem we get
$$(\hat{G},\phi) = \int _{\Bbb R^n} G(x) \hat{\phi}(x)\ dx =
\int_{\Bbb R^n} \Big(\int_\Omega |(x,\xi)|^q d\mu(\xi)\Big)
\hat{\phi}(x)\ dx =$$
$$\int_\Omega d\mu(\xi) \int_{\Bbb R^n} |(x,\xi)|^q \hat{\phi}(x)\ dx =
(2\pi)^{n-1} C_q \int_\Omega  d\mu(\xi)
\int_{\Bbb R} |t|^{-1-q} \phi(t\xi)\ dt. \qed$$
\enddemo

\bigbreak

\subheading{Remark 2} Lemma 2 was proved in \cite{14} in a 
slightly different form, and it was used there to give
a new Fourier transform proof of the following well-known
uniqueness theorem (see introduction for the history of 
the problem  and other applications): if $q>0,$ $q$ is not 
an even integer, and $\mu$ and $\nu$ are symmetric 
measures with bounded variation on $\Omega$ so that,
for every $x\in \Omega$                             
$$ \int_\Omega |(x,\xi)|^q\ d\mu(\xi) =
\int_\Omega |(x,\xi)|^q\ d\nu(\xi),\tag{5}$$
then $\mu=\nu.$ To see that, it is enough to apply
Lemma 2 to the test functions of the form 
$\phi(x)= u(t) v(\xi),$ where $x=t\xi,\ t>0,\ \xi\in \Omega,$
$u$ is any test function on $\Bbb R$ with $0\notin supp(u),$
and $v$ is any even infinitely differentiable function on
the sphere $\Omega.$ For such functions $\phi,$ we have
$\int_{\Bbb R} |t|^{-1-q} \phi(t\xi)\ dt = v(\xi)
\int_{\Bbb R} |t|^{-1-q} u(t)\ dt.$ Since the Fourier transforms 
of both sides of (5) are equal and have the property of Lemma 2,
we derive from (4) that
$\int_\Omega v(\xi) d\mu(\xi) = \int_\Omega v(\xi)\ d\nu(\xi)$
for any infinitely differentiable function $v$ on $\Omega,$ which
implies $\mu=\nu.$  Note that if $q$ is an even integer the 
uniqueness theorem fails to be true because only a finite number of
moments of the measures $\mu$ and $\nu$ must be equal.

\bigbreak

Now we are ready to show the connection between the Fourier
transform and the Blaschke-Levy representation.

\proclaim{Theorem 1} Let $H$ be a continuous, non-negative, even 
homogeneous function of degree 1 on $\Bbb R^n.$
Suppose that, for some $q>0$ which is not an even integer, 
$(H^q)^{\wedge}$ is a function on $\Bbb R^n\setminus\{0\}$ 
so that $(H^q)^{\wedge}\vert_{\Omega}$ belongs to the space 
$L_1(\Omega).$
Then the function $H^q$ admits the Blaschke-levy representation
with the exponent $q,$ and the corresponding function $b\in L_1(\Omega)$ 
is given by 
$b(\xi) = \big(1/(2(2\pi)^{n-1}C_q)\big)(H^q)^{\wedge}(\xi)$ 
for every $\xi\in \Omega.$
\endproclaim

\demo{Proof} 
The Fourier transform of an even homogeneous (of degree $q$)
function $H^q$ is an even homogeneous distribution of degree $-n-q.$
Fix any even test function $\phi$ with 
$0\notin supp(\phi).$ Since we know that $(H^q)^{\wedge}$ is
a function on $\Bbb R^n\setminus\{0\}$ whose restriction to
the sphere is an $L_1$-function, we can write the value of the
distribution $(H^q)^{\wedge}$ at the test function $\phi$ as
an integral, and then pass to the spherical coordinates: 
$$\big( (H^q)^{\wedge},\phi \big) = 
\int_{\Bbb R^n} (H^q)^{\wedge}(x)\ \phi(x)\ dx = 
(1/2)\int_\Omega (H^q)^{\wedge}(\xi)\ d\xi 
\int_{\Bbb R} |t|^{-1-q}\phi(t\xi)\ dt.\tag{6}$$

Put $b(\xi) = \big(1/(2(2\pi)^{n-1}C_q\big)(H^q)^{\wedge}(\xi)$ 
for every $\xi\in \Omega,$ and let us show that this function $b$
provides the 
equality (1). Since $b\in L_1(\Omega),$ the integral
$\int_\Omega |(x,\xi)|^q b(\xi)\ d\xi$ is a homogeneous function
(of the variable $x\in \Bbb R^n)$ of  degree $q$ whose restriction 
to the sphere is an $L_1$-function. By Lemma 2,
$$\Big(\big(\int_\Omega |(x,\xi)|^q b(\xi)\ d\xi \big)^{\wedge}, 
\phi \Big) =
(2\pi)^{n-1}C_q\int_\Omega b(\xi) \int_{\Bbb R} |t|^{-1-q} \phi(t\xi)\ dt.
\tag{7}$$
Because of the definition of the function $b,$ the right-hand sides
of (6) and (7) are equal. Since $\phi$ is an arbitrary  even test
function supported in $\Bbb R^n\setminus \{0\},$ the even functions
$(H^q)^{\wedge}$ and 
$\big(\int_\Omega |(x,\xi)|^q b(\xi)\ d\xi \big)^{\wedge}$ are 
equal distributions in $\Bbb R^n\setminus \{0\}.$ Therefore,
$H^q$ and $x\to \int_\Omega |(x,\xi)|^q b(\xi)\ d\xi$ are functions
in $\Bbb R^n$ which can differ by a polynomial only 
(see \cite{9, p. 119}). Since
both of those functions are even homogeneous of the order $q,$
and $q$ is not an even integer, we conclude that the polynomial 
must be equal to zero, and we have (1). The uniqueness 
follows from Remark 2. \qed \enddemo

\bigbreak

We end this section by showing that the Fourier transform 
of a homogeneous function of degree $p\le -n+1$ can be expressed
in terms of the function itself. We have to treat the cases 
$p<-n+1$ and $p=-n+1$ separately.

\proclaim{Lemma 3} Let $p<-n+1$ so that $-n-p$ is not an even  
integer, and let $f$ be an even homogeneous function of degree
$p$ on $\Bbb R^n\setminus \{0\},\ n>1$ such that 
$f\vert_\Omega \in L_1(\Omega).$ 
Then for every $\xi\in \Bbb R^n$
$$\hat{f}(\xi) = {{\pi}\over {C_{-n-p}}} 
\int_\Omega |(\theta,\xi)|^{-n-p} f(\theta)\ d\theta,\tag{8}$$
so $\hat{f}\vert_\Omega \in L_1(\Omega),$ and  
$\|\hat{f}\vert_\Omega\|_1 \le (\pi W_{-n-p}/C_{-n-p}) 
\|f\vert_\Omega\|_1.$
Also if $f\vert_\Omega \in L_{\infty}(\Omega)$ then
$\hat{f}\vert_\Omega \in L_\infty(\Omega)$ and
$\|\hat{f}\vert_\Omega\|_\infty 
\le (\pi W_{-n-p}/C_{-n-p}) \|f\vert_\Omega\|_\infty.$
\endproclaim

\demo{Proof} Since $f\vert_\Omega\in L_1(\Omega)$ 
and $-n-p > -1,$ Remark 1 implies that the right-hand side 
of (8) is a homogeneous 
function of degree $-n-p$ whose restriction to the sphere
is an $L_1$-function. 
Let $\phi$ be an even test function with
$0\notin supp(\hat{\phi}).$ Switching to the spherical coordinates 
and using the fact that $f$ is even homogeneous we get
$$(\hat{f},\phi)= \int_{\Bbb R^n} f(z) \hat{\phi}(z)\ dz=
(1/2)\int_\Omega \int_{\Bbb R} 
f(t\theta)\ |t|^{n-1} \hat{\phi}(t\theta)\ dt \ d\theta=$$ 
$$(1/2)\int_\Omega f(\theta)\ d\theta 
\int_{\Bbb R} |t|^{n+p-1} \hat{\phi}(t\theta)\ dt. 
\tag{9}$$
Now we apply Lemma 1 
with $q=-n-p.$ Recall that $(\hat\phi)^\wedge = (2\pi)^n \phi.$
The right-hand side of (9) is equal to
$${{(2\pi)^n}\over{2 (2\pi)^{n-1}C_{-n-p}}}
\int_\Omega f(\theta)\ d\theta 
\int_{\Bbb R^n} |(\theta,\xi)|^{-n-p} \phi(\xi)\ d\xi=$$
$${\pi\over{C_{-n-p}}}
\Big( \int_\Omega |(\theta,\xi)|^{-n-p} f(\theta)\ d\theta, \phi \Big).$$
Since $\phi$ is an arbitrary even test function 
with $0\notin supp(\hat{\phi})$ we conclude (similarly
to the end of the proof of Theorem 1) that the functions 
$\hat{f}(\xi)$ and $\xi \to 
(\pi/C_{-n-p})\int_S |(\theta,\xi)|^{-n-p} f(\theta)\ d\theta$
are even homogeneous functions of the order $-n-p$ which are 
equal up to an even homogeneous polynomial, and that polynomial must be 
equal to zero because the number $-n-p$ is not an even integer.
So we get (8), 
and the inequalities for the norms follow.\qed \enddemo

\bigbreak

\proclaim{Lemma 4} 
Let $f$ be an even homogeneous function of degree $-n+1$
on $\Bbb R^n\setminus \{0\},\ n>1$ so that 
$f\vert_\Omega\in L_1(\Omega).$ Then, for every $\xi\in \Omega,$
$$\hat{f}(\xi) = 
\pi \int_{\Omega\cap \{(\theta,\xi)=0\}} f(\theta)\ d\theta.$$
In particular, if $f\vert_\Omega\in L_\infty(\Omega)$ then 
$\hat{f}\vert_\Omega\in L_\infty(\Omega),$ and 
$$\|\hat{f}\vert_\Omega\|_\infty \le 
(2\pi^{(n+1)/2}/\Gamma((n-1)/2)\|f\vert_\Omega\|_\infty.$$
\endproclaim

\demo{Proof} Because of the connection between
the Fourier transform and the Radon transform,
for every even test function $\phi$ and every 
$\theta\in \Omega,$ the Fourier transform of the function
$t\to \hat{\phi}(t\theta)$ at zero is equal to 
$\int_{\Bbb R} \hat{\phi}(t\theta)\ dt = 
2\pi\int_{(\theta,\xi)=0} \phi(\xi)\ d\xi.$ Also
the Fourier transform of the $\delta$-function (defined by
$(\delta, \phi)=\phi(0)$) is the constant function 
$h(t)=1.$ 
Therefore, passing to the spherical coordinates we get
$$(\hat{f},\phi) = \int_{\Bbb R^n} f(x)\hat{\phi}(x)\ dx =
\int_\Omega \int_0^\infty 
f(t\xi) t^{n-1} \hat{\phi}(t\xi)\ dt\ d\xi =$$
$$(1/2)\int_\Omega f(\theta)\ d\theta 
\int_{\Bbb R} \hat{\phi}(t\theta)\ dt= 
\pi \int_\Omega f(\theta)\ d\theta
\int_{(\theta,\xi)=0} \phi(\xi)\ d\xi =$$
$$\pi \int_{\Bbb R^n} 
\Big(\int_{\Omega\cap \{(\theta,\xi)=0\}} f(\theta)\ d\theta\Big)
\phi(\xi)\ d\xi,$$ 
and the result follows since $\phi$ is an arbitrary
even test function. \qed \enddemo

\head 3. The inverse formula.  \endhead

Theorem 1 gives a condition for the existence of 
the Blaschke-Levy representation and the inverse formula
in terms of the Fourier transform of the original
function. Though this criterion has a few applications
(see \cite{14}), it is often difficult to calculate the 
Fourier transform. However, using Lemmas 3 and 4 we 
can replace the Fourier transform condition by a
condition in terms of the derivatives of the original
function which is sometimes more convenient for applications.
 
Let us explain what is going to happen. Suppose we want to find
the Blaschke-Levy representation for a function $H.$ Theorem 1 reduces
this problem to calculating the Fourier transform of $H.$
Instead of doing that, let us consider the distribution $\Delta^k H,$
where $\Delta$ is the Laplace operator and $k$ is an integer so that the
distribution $\Delta^k H$ is homogeneous of degree less or equal
than $-n+1.$ The Fourier transform of $\Delta^k H$ has (up to a sign)
the same restriction to the sphere as the Fourier transform of $H.$
On the other hand, by Lemma 3 (or Lemma 4) if $\Delta^k H$ is an 
$L_1$-function on the sphere
so is its Fourier transform, and there is a simple formula 
expressing the Fourier transform of $\Delta^k H$ in terms of the 
function itself. That is why we can replace the Fourier transform 
condition for the existence of the Blaschke-Levy representation by a
condition in terms of the function $\Delta^k H.$
\bigbreak

First, let us consider the case where $q$ is not an integer.

\proclaim{Theorem 2} Let $q>0,$ $q$ is not an integer, and
let  $H$ be a continuous, non-negative, 
even homogeneous function of degree 1 on $\Bbb R^n,\ n>1.$
Suppose that $\Delta^k H^q$ is a function in 
$\Bbb R^n\setminus\{0\}$ so that 
$(\Delta^k H^q)\vert_\Omega \in L_1(\Omega),$ 
where $k=(n+[q])/2$ if $n+[q]$ is an even integer, $k=(n+[q]+1)/2$
if $n+[q]$ is an odd integer, and differentiation is considered 
in the sense of distributions.
Then the function $H^q$ admits the Blaschke-Levy
representation (1) with the exponent $q,$ where the function 
$b\in L_1(\Omega)$ can be calculated by
$$b(\xi) = (-1)^k {\pi\over{2(2\pi)^{n-1}C_{-n-q+2k}C_q}} 
\int_\Omega |(\theta,\xi)|^{-n-q+2k} 
(\Delta^k H^q)(\theta)\ d\theta$$
for every $\xi\in \Omega.$ 
Moreover,
$$\|b\|_1 \le  {{\pi W_{-n-q+2k}}\over {2(2\pi)^{n-1} C_{-n-q+2k}C_q}}  
\|(\Delta^k H^q)\vert_\Omega\|_1.$$ 
If the function $(\Delta^k H^q)\vert_\Omega$ belongs to 
$L_\infty(\Omega)$ then $b\in L_\infty(\Omega)$ and
$$\|b\|_\infty \le {{\pi W_{-n-q+2k}}\over{2(2\pi)^{n-1}
C_{-n-q+2k}C_q}}  \|(\Delta^k H^q)\vert_\Omega\|_\infty.$$ 
\endproclaim

\demo{Proof} Since the function $H^q$ is even homogeneous, 
the distribution $\Delta^k H^q$ is even homogeneous of the order
$q-2k<-n+1.$ Also $-n-q+2k$ is not an even integer, so  $\Delta^k H^q$
satisfies the conditions of Lemma 3. By Lemma 3, the Fourier transform
$$(\Delta^k H^q)^{\wedge} (\xi) = 
{\pi\over{C_{-n-q+2k}}} 
\int_\Omega  |(\theta,\xi)|^{-n-q+2k} (\Delta^k H^q)(\theta)
\ d\theta$$
is an $L_1$-function on $\Omega.$ Because of the connection between
the Fourier transform and differentiation we have
$$(\Delta^k H^q)^{\wedge} (\xi) = 
(-1)^k (\xi_1^2+...+\xi_n^2)^k (H^q)^{\wedge}(\xi),$$
and the restrictions to the sphere of $(\Delta^k H^q)^{\wedge}$
and $(-1)^k (H^q)^{\wedge}$ are equal. In particular,
$(H^q)^{\wedge}\vert_\Omega \in L_1(\Omega).$
This means that we can apply Theorem 1, and the result follows.
\qed \enddemo

\bigbreak

If $q$ is an even integer the uniqueness in the Blaschke 
representation fails (as mentioned in Remark 2). Therefore,
it remains to consider the case where $q$ is an odd integer.

First, suppose that the dimension $n$ is even. Then  we apply 
the Laplace operator to the function $H^q$ until it becomes 
a homogeneous function
of degree $-n+1,$ and then we use Lemma 4 instead of Lemma 3. 
The rest of the proof of Theorem 3 is similar to that of
Theorem 2. 

\proclaim{Theorem 3} Let $n\in \Bbb N$ be an even integer,
$q>0$ be an odd integer, and $H$ is a continuos, non-negative, 
even homogeneous
function of degree 1 on $\Bbb R^n,\ n>1$.
Suppose that $\Delta^{(n+q-1)/2} H^q$ is a function in 
$\Bbb R^n\setminus\{0\}$ so that 
$(\Delta^{(n+q-1)/2} H^q)\vert_\Omega \in L_1(\Omega),$ 
where differentiation is considered in the sense of distributions.
Then the function $H^q$ admits the Blaschke-Levy representation (1)
with the exponent $q,$
and the corresponding function $b\in L_1(\Omega)$ is given by
$$b(\xi) = {{(-1)^{(n+q-1)/2}\pi}\over{(2\pi)^{n-1}C_q }}
\int_{\Omega\cap \{(\theta,\xi)=0\}} 
\Delta^{(n+q-1)/2} H^q(\theta)\ d\theta$$
for every $\xi\in \Omega.$ Moreover, if $b\in L_\infty(\Omega)$ then
$$\|b\|_\infty \le {{2\pi^{(n+1)/2}}\over {\Gamma((n-1)/2)(2\pi)^{n-1}C_q }}
\|(\Delta^{(n+q-1)/2} H^q)\vert_\Omega\|_\infty.$$ 
\endproclaim

\bigbreak

In the case where $q$ and $n$ are both odd integers, the 
technique of this paper does not work. The reason is that the
polynomials, which appear at the end of the proofs of Theorem 1
and Lemma 3 (and can easily be eliminated in those cases),
start playing active role when $n+q$ is an even integer.
To illustrate this, let us just note that, for the Euclidean
norm $\|x\|_2$ in $\Bbb R^n$ with $n$ being an odd integer, 
the distribution
$\Delta^2 \|x\|_2$ vanishes everywhere in $\Bbb R^n \setminus \{0\},$
and, therefore, it is a linear combination of the derivatives of the 
$\delta$-function. So in the case where $q$ and $n$ are odd integers,
the Fourier transform of $\Delta^k H^q$ not always can be expressed
in terms of the restriction of the function $\Delta^k H^q$ to the
sphere.

\bigbreak

Let us give a scheme of how Theorem 3 works in the case 
where $q=1,$ $n$ is 
an even integer, and the function $H$ is of the form 
$H(x)= P_m(x)\|x\|_2^{-m+1}$ on $\Bbb R^n,$
where $P_m$ is an even homogeneous polynomial of degree 
$m > 0.$ 
First, by Euler's formula for homogeneous
functions, we have $\sum x_i (\partial P_m/\partial x_i) = mP_m,$
and, for every $\beta,$ 
$$\Delta (P_m\|x\|_2^\beta) = \Delta (P_m) \|x\|_2^{\beta} + 
\beta (n+2m+\beta-2) P_m \|x\|^{\beta-2}.$$
Iterating the latter formula one can calculate 
$\Delta^k H$ for every k,
and find the polynomial which is the restriction 
of $\Delta^k H$ to the sphere. Now the problem of finding
the Blaschke-Levy representation for the function $H$ is
reduced to calculating the integrals of the form
$$\int_{\Omega\cap \{(x,\xi)=0\}} 
x_1^{\alpha_1}\dots x_n^{\alpha_n}\ dx,\tag{10}$$
where $\alpha_i$ are even integers and $\xi\in \Omega.$ 
To calculate these integrals we use an argument
similar to that of \cite{17}.
Namely, we start with the equality
$$\|\xi\|_2^{\alpha_1+\dots+\alpha_n - 1} = 
(1/W_{\alpha_1+\dots+\alpha_n - 1}) 
\int_\Omega |(x,\xi)|^{\alpha_1+\dots+\alpha_n-1}\ dx.$$
Differentiating this equality we see that the integral (10) is equal to
$${{\partial^{\alpha_1+\dots+\alpha_n}
\|\xi\|_2^{(\alpha_1+\dots+\alpha_n-1)/2}}\over
{\partial \xi_1^{\alpha_1} \dots \partial \xi_n^{\alpha_n}}}
{{W_{\alpha_1+\dots+\alpha_n - 1}}\over {(\alpha_1+\dots+\alpha_n - 1)!}},$$
where the derivative is calculated at the point $\xi\in \Omega.$
In Section 4 we give a numerical example.

\medbreak

This calculation includes differentiation
only. A different way of calculating the function
$b$ is to find the spherical harmonics expansion of the 
polynomial $P_m,$ and then use Rodriguez's formula (see \cite{22}
for the properties of spherical harmonics).

If $H$ is not of a polynomial form, it is sometimes impossible
to calculate $b$ precisely using our inverse formulae. 
However, Theorems 2 and 3 give
estimates for the $L_1$ and $L_\infty$- norms of the function
$b$ with computable constants. This seems to be an advantage 
of our approach over the one using spherical harmonics where
the constants appear as the sums of certain series'.

\bigbreak

\head 4. A characterization of subspaces of $L_q.$ \endhead

The question of how to check whether a given space is isometric 
to a subspace of $L_q$ is a matter of an old problem raised by 
P.Levy \cite{19}. 
In \cite{19} P.Levy showed that an $n$-dimensional 
space is isometric to a subspace of $L_q$ if and only if its norm 
admits the Blaschke-Levy representation with the exponent $q$
(and with a non-negative
measure in place of the function $b.)$ Bretagnolle, Dacunha-Castelle
and Krivine \cite{5} proved that,
for $0<q\le 2,$ a Banach space is isometric to a subspace of $L_q$
if and only if the function $\exp(-\|x\|^q)$ is positive definite,
and, in particular, showed that the space $L_p$ embeds isometrically
into $L_q$ if $0<q<p\le 2.$ Another criterion involving the Fourier 
transform (which,
in fact, is our Theorem 1 in a slightly stronger form) was given
in \cite{14}, \cite{15}: for any 
$q\in (0,\infty)\setminus \{even\ integers\},$
an $n$-dimensional space is isometric to a subspace of $L_q$ if and
only if the restriction of the Fourier transform of $\|x\|^q \Gamma(-q/2)$ 
to the sphere $\Omega$ is a finite Borel (non-negative)
measure on $\Omega.$
Though the Fourier transform criteria work for certain spaces,
calculating the Fourier transform of a norm precisely is not
always possible. That is why
a condition involving the derivatives of the norm instead of the 
Fourier transform could be useful.
A necessary condition in terms of the derivatives
of the norm was given by Zastanvy \cite{29} who proved
that a three dimensional space 
is not isometric to a subspace of $L_q$ with $0<q\le 2$ if there exists 
a basis $e_1,e_2,e_3$ so that the function
$$(y,z)\mapsto \|xe_1 + ye_2 + ze_3\|^{'}_x (1,y,z)/\|e_1+ye_2+ze_3\|,
\ y,z \in \Bbb R$$
belongs to the space $L_1(\Bbb R^2).$ 

In this section, we use the inverse formula for the 
Blaschke-Levy representation to give a sufficient condition for 
the existence of isometric embedding of a space into $L_q$ 
which is formulated in terms of the Laplace operator of the
norm.

\bigbreak

We start with a well-known fact which explains the connection
between the Blaschke-Levy representation and isometric embedding
into $L_q.$  

\proclaim{Lemma 5} Let $q$ be a positive number which is 
not an even integer, $(X, \|\cdot\|)$ be an $n$-dimensional 
space, and suppose that the function $\|x\|^q$ admits the 
Blaschke-Levy representation with a function $b\in L_1(\Omega):$  
for every $x\in \Bbb R^{n}$,
$$ \|x\|^{q} = \int_{\Omega} |( x,\xi )|^{q}\ b(\xi)\ d\xi. \tag{11}$$
Then $X$ is isometric to a subspace of $L_q$ if and only if
$b$ is a non-negative (not identically zero) function.
\endproclaim 

\demo{Proof} If $b$ is a non-negative function we can assume without 
loss of generality that $\int_{\Omega} b(\xi)\ d\xi\ =\ 1.$ Choose any
measurable (with respect to Lebesgue measure) functions $f_1,\dots,f_n$
on $[0,1]$ so that their joint distribution is the measure $b(\xi)d\xi$ 
on the sphere $\Omega.$ Then, by (11),
the operator $x\mapsto \sum x_i f_i,\ x\in \Bbb R^n$ is an isometry 
from $X$ to $L_q([0,1]).$

Conversely, if $X$ is a subspace of $L_q([0,1])$ choose any functions 
$f_{1},...,f_{n}\in L_q$ which form a basis in $X,$ and let $\mu$ be 
the joint distribution of the functions $f_{1},...,f_{n}$ with 
respect to Lebesgue measure. Then, for every $x\in R^{n}$,
$$\|x\|^{q}=\|\sum_{k=1}^{n} x_{k}f_{k}\|^{q}=\int_{0}^{1} 
|\sum_{k=1}^{n}x_{k}f_{k}(t)|^{q} dt=$$
$$\int_{R^{n}} |( x,\xi)|^{q}\ d\mu (\xi)=\int_{\Omega} 
|(x,\xi)|^{q}\ d\mu_\Omega (\xi)\tag{12}$$
where $\mu_\Omega$ is the projection of $\mu$ to the sphere.(For every 
Borel subset $A$ of $\Omega$, 
$\mu_\Omega(A)=(1/2)\int_{\{tA,t\in R\}} \|x\|_{2}^{q} d\mu(x) ).$
It follows from (11) and (12) that
$$\int_{\Omega} |( x,\xi )|^{q}\ b(\xi)\ d\xi = 
\int_{\Omega} |( x,\xi )|^{q}\  d\mu_\Omega(\xi)$$
for every $x\in \Bbb R^n.$ Since $q$ is not an even integer,
we can apply the uniqueness theorem for measures on the sphere
(see Remark 2) to show that $ d\mu_\Omega(\xi) = b(\xi)\ d\xi$ which
means that $ b(\xi)\ d\xi$ is a measure, and the function $b$
is non-negative.   
\qed \enddemo

\bigbreak
In view of Lemma 5, the inverse formulae from Section 3
lead to the following criteria of isometric embedding into $L_q.$
\medbreak
First, if $q$ is not an integer
we use  Lemma 5 and Theorem 2: under the assumption that  
$(\Delta^k \|x\|^q)\vert_\Omega \in L_1(\Omega),$
an $n$-dimensional normed (quasi-normed)
space $(\Bbb R^n,\|\cdot\|)$ embeds isometrically in $L_q$ 
if and only if 
$$\xi\to {{(-1)^k}\over{C_{-n-q+2k}C_q}} 
\int_\Omega |(\theta,\xi)|^{-n-q+2k} 
(\Delta^k \|\theta\|^q)\ d\theta,$$
is a non-negative function on $\Omega,$ where $k$ is as in Theorem 2.
If for some reason it is impossible to calculate the latter
integral precisely, one can use the following sufficient condition:
if the function $((-1)^k/(C_{-n-q+2k}C_q)) 
\Delta^k \|x\|^q$ is non-negative on $\Omega$ and its restriction to 
$\Omega$ belongs to $L_1(\Omega)$ then 
the space $(\Bbb R^n,\|\cdot\|)$ embeds isometrically into $L_q.$

\medbreak

If $q$ is an odd integer and the dimension $n$ is an even integer,
similar criteria follow from Lemma 5 and Theorem 3. Under the assumption
that $(\Delta^{(n+q-1)/2} \|x\|^q)\vert_\Omega \in L_1(\Omega),$
a space $(\Bbb R^n,\|\cdot\|)$ is isometric
to a subspace of $L_q$ if and only if 
$$\xi\to {{(-1)^{(n+q-1)/2}}\over{C_q }}
\int_{\Omega\cap (\theta,\xi)=0} 
\Delta^{(n+q-1)/2} \|\theta\|^q\ d\theta$$
is a non-negative function on $\Omega.$
The related sufficient condition is that
the function $((-1)^{(n+q-1)/2}/C_q) \Delta^{(n+q-1)/2} \|x\|^q$
is a non-negative $L_1$-function on $\Omega.$

\subheading{Example 1} Consider the 
function  $\|x\| = \|x\|_2 + \lambda x_1^2\|x\|_2^{-1}$
which is an even homogeneous function of degree 1 on $\Bbb R^n.$
For which values of $\lambda$ does the space $(\Bbb R^4,\|\cdot\|)$
embed isometrically in $L_1 ?$ An equivalent question asks 
for the values of $\lambda$ for which the polar set to
$\{x: \|x\|\le 1\}$ is a zonoid (see \cite{4} for the connection
between zonoids and embedding into $L_1.$)

Let us apply Theorem 3 with $q=1,\ n=4$ to find the function
$b$ corresponding to $H(x)=\|x\|.$ Since $(n+q-1)/2 =2$
we calculate                               
$$\Delta^2 \|x\| = -3\|x\|_2^{-3} + \lambda (-12\|x\|_2^{-1}+
45x_1^2 \|x\|_2^{-3}).$$
Therefore, $\Delta^2 H\vert_\Omega = -3-12\lambda + 45\lambda x_1^2.$
Also $C_1=-1,$ and, by Theorem 3, for every $\xi\in \Omega$
$$b(\xi) = {1\over{8\pi^2}}\int_{\Omega\cap \{(\theta,\xi)=0\}}
(3+12\lambda-45\lambda x_1^2)\ dx.$$
To calculate the integral note that 
$\int_{\Omega\cap \{(\theta,\xi)=0\}} x_1^2\ dx$ is equal to
the second derivative by $\xi_1$ of the integral
$\int_\Omega |(x,\xi)|\ dx = W_1 \|\xi\|_2.$ 
Also the surface area of the 3-dimensional sphere 
${\Omega\cap \{(\theta,\xi)=0\}}$ is equal to $2\pi^{3/2}/\Gamma(3/2).$ 

Finally, $b(\xi) = (1/(8\pi))(4 - 8\lambda + 24 \lambda \xi_1^2).$
Clearly, $b$ is a non-negative function if and only if
$-1/4 \le \lambda \le 1/2,$ and these are all the values of $\lambda$
for which the space embeds in $L_1.$

\bigbreak

\subheading{Example 2} Let $\|x\| = \|x\|_2 + \lambda P(x),$
where $P$ is an even homogeneous function of degree 1 on
$\Bbb R^n,\ n$ is an even integer, 
and $P\vert_\Omega \in C^{n/2}(\Omega).$
To find the values of
$\lambda$ for which $(\Bbb R^n,\|\cdot\|)$ embeds 
isometrically in $L_1,$ we calculate 
$$\Delta^{n/2}\|x\| = (-1)^{(n-2)/2}(n-1)!! (n-3)!!\|x\|_2^{-n+1} +
\lambda \Delta^{n/2}P .$$
The sufficient condition formulated above shows that 
the space embeds in $L_1$ if 
$(n-1)!! (n-3)!!- (-1)^{(n-2)/2}\lambda (\Delta^{n/2}P)\vert_\Omega$ 
is a non-negative function (note that $C_1 = -1.$) Hence, if
$$|\lambda|\le {{(n-3)!!(n-1)!!}
\over{\|(\Delta^{n/2}P)\vert_\Omega\|_\infty}}.$$ 
then the space embeds in $L_1.$
\bigbreak

\head 5. A Fourier transform formula for the central
sections of star bodies \endhead

Let $K$ be a centrally symmetric star body in $\Bbb R^n$ so that   
the norming functional  $\|x\| = \min\{a>0 : x\in aK\},
\ x\in \Bbb R^n$ generated by $K$ is a continuous,  
non-negative, even homogeneous function of degree 1
on $\Bbb R^n.$ It is easy to
see that, for every $\xi$ in the unit sphere $\Omega,$ 
the $(n-1)$-dimensional volume of the section of $K$ by 
the hyperplane $\xi^{\bot}=\{(x,\xi)=0\}$ satisfies the equality
$${{{Vol}_{n-1}(K\cap \xi^{\bot})}\over {Vol_{n-1}(B_{n-1})}} =
{{\int_{\Omega\cap \xi^{\bot}} \|x\|^{-n+1}\ dx}\over{A_{n-1}}},
\tag{13}$$
where $Vol_{n-1}(B_{n-1})= \pi^{(n-1)/2}/\Gamma((n+1)/2)$ is the
volume of the Euclidean unit ball $B_{n-1}$ in $\Bbb R^{n-1},$
and $A_{n-1}= 2\pi^{(n-1)/2}/\Gamma((n-1)/2))$ is the surface area
of the Euclidean unit sphere in $\Bbb R^{n-1}.$

The integral in the right-hand side of (13)
is equal to the integral in Lemma 4 with
$f(x)=\|x\|^{-n+1}.$ Therefore, Lemma 4 and (13) imply the 
following Fourier transform formula for the volume
of central sections of $K:$

\proclaim{Theorem 4} For every $\xi\in \Omega,$
$${Vol}_{n-1}(K\cap \xi^{\bot}) = {1\over{\pi (n-1)}}
(\|x\|^{-n+1})^{\wedge}(\xi).$$
\endproclaim

\bigbreak

The Fourier transforms of powers of different norms 
have been calculated in \cite{18} (for the $\ell_\infty^n$-norm),
\cite{15} (for the $\ell_p^n$-norms), \cite{6} (for the Lorentz
norm). In view of Theorem 4, one can use those calculations 
to obtain formulae for the volume of central sections.
For example, the Fourier transform of the functions
of the form $f(\|x\|_\infty)$ was calculated in \cite{18},
where $\|x\|_\infty$ stands for the norm of 
the space $\ell_\infty^n,$ and $f$ belongs to a large class
of functions on $\Bbb R.$ (Note
that a multiplier $(-1)^{n-1}$ is missing in the formula
in \cite{18}.)
If we apply the formula from \cite{18} to the function
$f(t)=|t|^p$ with $p\in (-1,0),$ use the formulae  
for the Fourier transform of the functions $|t|^p$ and 
$|t|^p sgn(t)$ (see \cite{8, p.173}), and then 
use analytic extension by $p,$ we get an expression for
the Fourier transform of $\|x\|_\infty^{-n+1}:$
for every $\xi \in \Bbb R^n$ with non-zero coordinates,
if the dimension $n$ is odd we have
$$(\|x\|^{-n+1})^{\wedge} (\xi) = $$
$${{(-1)^{(n-1)/2} 2^{-n+1} \sqrt{\pi}\ \Gamma((-n+2)/2)}\over
{\Gamma((n-1)/2)\prod_{k=1}^n \xi_k}}
\sum_\delta \delta_1\dots\delta_n 
\big|\sum_{j=1}^n \delta_j\xi_j\big|^{n-1}
sgn(\sum_{j=1}^n \delta_j\xi_j).$$
If the dimension $n$ is even we have
$$(\|x\|^{-n+1})^{\wedge} (\xi) =
{{(-1)^{(n-2)/2} 2^{-n+1} \sqrt{\pi}\ \Gamma((-n+3)/2)}\over
{\Gamma(n/2)\prod_{k=1}^n \xi_k}}
\sum_\delta \delta_1\dots\delta_n 
\big|\sum_{j=1}^n \delta_j\xi_j\big|^{n-1}.$$
The outer sum is taken over all changes of sign 
$\delta = (\delta_1,\dots,\delta_n),\ \delta_j = \pm 1.$
These formulae, in conjunction with Theorem 4, imply simple
formulae for the volume of central sections of the cube
$[-1,1]^n.$ Previously, similar formulae were obtained using 
probabilistic arguments specifically designed for the cube.
Ball \cite{2} has shown that the exact lower and upper
bounds for the volume of central sections
of the unit ball of the space $\ell_\infty^n$ are
$2^n$ and $2^n\sqrt{2},$ respectively. We refer the reader
to \cite{7} for a historical survey and more information
about sections.

\bigbreak 

Meyer and Pajor \cite{21} have proved that the minimal section
of the unit ball of the space $\ell_1^n$ is the one 
perpendicular to the vector $(1,1,...,1),$ and
the maximal section is perpendicular to the vector
$(1,0,...,0).$ They also showed that, for the unit balls of 
the spaces $\ell_p^n$ with $1<p<2,$ the upper bound 
occurs in the same direction as for $p=1,$ and 
raised the question of whether 
the situation is the same for the lower bound.

We end this paper by confirming the conjecture of Meyer
and Pajor. First, let us compute the Fourier transform 
of the functions $\|x\|_p^\beta,$ where $\|x\|_p$ stands 
for the norm of the space $\ell_p^n.$
Denote by $\gamma_p$ the Fourier transform of
the function $z\to \exp(-|z|^p),\ z\in \Bbb R.$
For $0<p\le 2,$ $\gamma_p$ is (up to a constant)
the density of the standard $p$-stable measure
on $\Bbb R,$ so $\gamma_p$ is a non-negative function.
For every $p>0,$
$$
\lim\limits_{t\to\infty} t^{1+p}\gamma_p(t)=2\Gamma (p+1) \sin (\pi p/2),
$$
so $\gamma_p$ decreases at infinity as $|t|^{-1-p}$ (see
\cite{30}). 
Also simple calculations show that 
$\gamma_p(0) = 2\Gamma(1+ 1/p),$ and  
$\int_0^\infty \gamma_p(t)\ dt = \pi.$
The following calculation is taken from \cite{15}.

\proclaim{Lemma 6}
Let $p>0$, $n\in\Bbb N$, $-n<\beta<pn$, $\beta/p \not\in \Bbb N\cup \{0\}$,
$\xi=(\xi_1, \ldots, \xi_n)\in \Bbb R^n$, $\xi_k\not= 0$, $1\leq k\leq n$.
Then
$$ (\|x\|_p^\beta)^{\wedge}(\xi) =
((|x_1|^p+\cdots +|x_n|^p)^{\beta/p})^\wedge
(\xi)=\frac{p}{\Gamma(-\beta/p)} 
\int_0^\infty t^{n+\beta-1}\prod_{k=1}^n
\gamma_p (t\xi_k)\, dt.
$$
\endproclaim

\demo{Proof}
Assume that $-1<\beta<0$.  By the definition of the $\Gamma$--function
$$
(|x_1|^p+\cdots +|x_n|^p)^{\beta /p}
=\frac{p}{\Gamma(-\beta/p)}\int_0^\infty y^{-1-\beta} \exp
(-y^p(|x_1|^p+\cdots + |x_n|^p))\, dy.
$$
For every fixed $y>0$, the Fourier transform of the function
$x\to \exp (-y^p(|x_1|^p+\cdots +|x_n|^p))$ at any point 
$\xi\in \Bbb R^n$ is equal to $y^{-n}\prod_{k=1}^n
\gamma_p (\xi_k/y).$ Making the change of variables  $t=1/y$  
we get
$$
((|x_1|^p+\cdots+|x_n|^p)^{\beta/p})^\wedge
(\xi)=\frac{p}{\Gamma(-\beta/p)}\int_0^\infty y^{-n-\beta-1} \prod_{k=1}^n
\gamma_p (\xi_k /y)\, dy=$$
$$\quad = \frac{p}{\Gamma(-\beta/ p)}\int_0^\infty t^{n+\beta-1}
\prod_{k=1}^n \gamma_p (t\xi_k)\, dt.\tag{14}$$

The latter integral converges if $-n< \beta <pn$ since the function
$t\to \prod_{k=1}^n \gamma_p (t\xi_k)$ decreases at infinity like
$t^{-n-np}$ (recall that $\xi_k\not= 0$, $1\leq k\leq n.)$ 

If $\beta$ is allowed to assume complex values then the both sides of
(14) are analytic functions of $\beta$ in the domain $\{-n<Re
\beta <np,\ \beta/p \not\in \Bbb N\cup \{0\}\}$.  These two functions admit
unique analytic continuation from the interval $(-1,0)$.  Thus the
equality (14) remains valid for all $\beta\in (-n,pn)$,
$\beta/p\not\in \Bbb N\cup \{0\}$ (see \cite{8} for
details of analytic continuation in such situations). \qed
\enddemo

\bigbreak

Now we can use Lemma 6 with $\beta=-n+1$ and Theorem 4
to get an expression for the volume of central sections.
Note that the condition of Lemma 6  that $\xi$ has non-zero
coordinates may be removed in Corollary 1 because the volume
of a section is a continuous function of $\xi.$
Denote by $B_p$ the unit ball of the space $\ell_p^n,\ p>0,\ n>1.$

\proclaim{Corollary 1} For every $p>0$ and $\xi\in \Omega,$
$${Vol}_{n-1}(B_p\cap \xi^\bot ) = {p\over{\pi(n-1)\Gamma((n-1)/p)}}
\int_0^\infty \prod_{k=1}^n
\gamma_p (t\xi_k)\, dt.$$
\endproclaim

For $p\in (1,2),$ the latter equality was established by Meyer and Pajor \cite{21} using a probabilistic argument.
Note that when $p\to \infty$ the formula (14) turns into the
expression used by Ball \cite{2} for the slices of the unit cube. 

\bigbreak

The following fact is a property of the functions
$\gamma_p$ with $p\in (0,2)$ only.

\proclaim{Lemma 7}  For every $p\in (0,2),$ the function
$\gamma_p(\sqrt{t})$ is log-convex on $(0,\infty).$
In other words, the function $\gamma_p^{'}(t)/(t\gamma_p(t))$
is increasing on $(0,\infty).$ Also, for every 
$k,m\in \Bbb N,\  k<m$ and every $t>0,$  we have
$\gamma_p^k(k^{-1/2}t)\gamma_p^{m-k}(0)
\ge  \gamma_p^m(m^{-1/2}t).$
\endproclaim

\demo{Proof} A well-known fact is that there exists a measure
$\mu$ on $[0,\infty)$ whose Laplace transform is equal to
$\exp(-t^{p/2}).$ This is a stable measure, and its properties
and asymptotic behavior of its density (which decreases
at infinity as $|t|^{-1-p/2},$ up to a constant) are described, 
for example, in \cite{30}. For every $z\in \Bbb R,$
we have
$$\exp(-|z|^p) = \int_0^\infty \exp(-uz^2)\ d\mu(u).$$
Calculating the Fourier transform of both sides of the
latter equality as functions of the variable $z,$
we get, for every $t\in \Bbb R,$
$$\gamma_p(t) = \sqrt{2\pi} 
\int_0^\infty u^{-1/2} \exp({{-t^2}\over{4u}})\ d\mu(u).$$
where the integral converges because of the asymptotics 
of the density of $\mu$ at infinity, as mentioned above.
Now the fact that 
$\gamma_p^2(\sqrt{(t_1+t_2)/2})\le \gamma_p(\sqrt{t_1})
\gamma_p(\sqrt{t_2})$
follows from the Cauchy-Schwarz inequality applied to the 
functions $\exp(-t_1/(8u))$ and $\exp(-t_2/(8u))$ and the measure
$u^{-1/2}\ d\mu(u),$ where $t_1,t_2$ are arbitrary positive numbers.
Therefore, the function $\gamma_p(\sqrt{t})$ is log-convex 
which implies the other two statements of Lemma 7. 
\qed \enddemo

\proclaim{Theorem 5} For every $p\in (0,2)$ and every $\xi\in \Omega,$
$${{p}\over{\pi(n-1)\Gamma((n-1)/p)}}
\int_0^\infty \gamma_p^n(t/\sqrt{n})\ dt \le        
{Vol}_{n-1}(B_p\cap \xi^\bot) \le $$
$${p\over{\pi(n-1)\Gamma((n-1)/p)}} 
\gamma_p^{n-1}(0) \int_0^\infty \gamma_p(t)\ dt=
{{2^{n-1} p (\Gamma(1+1/p))^{n-1}}\over{(n-1)\Gamma((n-1)/p)}}$$
with the left inequality turning into an equality 
if and only if $|\xi_i|=1/\sqrt{n}$ for every $i,$
and the upper bound occurs if and only if one of the coordinates 
of the vector $\xi$ is equal to $\pm 1$ 
and the others are equal to zero.
\endproclaim

\demo{Proof} Consider the function
$$F(\xi_1,\dots,\xi_n) = 
\int_0^\infty \gamma_p(t\xi_1)\cdots\gamma_p(\xi_n)\ dt +
\lambda (\xi_1^2+\dots+\xi_n^2 - 1),$$
where $\lambda$ is the Lagrange multiplier.
It suffices to find the maximal and minimal value
of the function $F$ in the positive octant under the
condition $\xi_2+\dots+\xi_n^2 = 1.$
To find the critical points of the function $F,$ 
we have to solve the system of equations
$${{\partial F}\over{\partial \xi_i}}(\xi) =
\int_0^\infty t\gamma_p^{'}(t\xi_i)\prod_{k\not= i}
\gamma_p (t\xi_k)\, dt + 2\lambda \xi_i = 0, $$
where $i=1,\dots,n.$ For each $i$ with $\xi_i\not= 0$ we can write 
the latter equality in the following form:
$$\int_0^\infty t{{\gamma_p^{'}(t\xi_i)}\over 
{\xi_i \gamma_p(t\xi_i)}}\prod_{k=1}^n
\gamma_p (t\xi_k)\, dt = - 2\lambda \tag{14}.$$
Since (by Lemma 7) the function 
$\gamma_p^{'}(t\xi_i)/(\xi_i \gamma_p(t\xi_i))$ is 
increasing and $\gamma_p$ is non-negative, we can  
have (14) for different values of $i$ simultaneously
only if the corresponding coordinates of the vector $\xi$
are equal.
Therefore, the critical points of the function $F$
are only those points $\xi$ for which some of the coordinates 
are zero, and the absolute values of the rest are equal. 
Hence, 
the problem is reduced to comparing the values of $F$
at the points $\xi^{(k)},\ k=1,...,n,$ where the first
$k$ coordinates of $\xi^{(k)}$ are equal to $1/\sqrt{k}$
and the last $n-k$ coordinates are equal to zero.
It follows from the inequality of Lemma 7 that the maximal 
value of $F$ on the sphere $\Omega$ occurs at the point $\xi^{(1)},$
and the minimal value is at the point $\xi^{(n)}.$ Now the
result of Theorem 5 follows from Corollary 1.
\qed \enddemo

\enddocument